\documentclass[12pt]{amsart}

\usepackage{amssymb}
\usepackage[latin1]{inputenc}

\title[Bulky Parabolic Subgroups]
{Special involutions and bulky parabolic subgroups in finite
  Coxeter groups}
\author[G.\ Pfeiffer and G.\ R\"{o}hrle]
{G\"otz Pfeiffer\ \  and\ \  Gerhard R\"{o}hrle}

\address%[G.\ Pfeiffer]
{Department of Mathematics, National University of Ireland, Galway, Ireland}
\email{goetz.pfeiffer@nuigalway.ie}
\address%[G.~R\"{o}hrle]
{School of Mathematics and Statistics, University of Birmingham, 
Birmingham B15 2TT, United Kingdom}
\email{ger@for.mat.bham.ac.uk}


\numberwithin{equation}{section}
\newtheorem{theorem}[equation]{Theorem}
\newtheorem{lemma}[equation]{Lemma}
\newtheorem{proposition}[equation]{Proposition}

\theoremstyle{definition}
\newtheorem{rem}[equation]{Remark}

\newcommand{\Size}[1]{\left|#1\right|}
\newcommand{\R}{\mathbb{R}}
\newcommand{\Z}{\mathbb{Z}}
\newcommand{\M}{\mathcal{M}}

\let\emptyset\varnothing

\begin{document}
\maketitle

\section{Introduction}
\label{sec:intro}

In \cite{FelderVeselov2004} Felder and Veselov 
considered the standard and twisted %(by complex conjugation) 
actions of a finite Coxeter group $W$  on the 
cohomology $H^*(\M_W)$ of the complement of the complexified hyperplane  
arrangement $\M_W$ of $W$.
The twisted action is obtained by combining the standard action with 
complex conjugation;
we refer the reader to  \cite{FelderVeselov2004} for precise
statements.
In a case by case argument, Felder and Veselov obtain a formula
for all  Coxeter groups $W$ for the standard action
\[
H^*(\M_W) \cong \sum_{\sigma \in X_W} 
(2 \cdot 1_{\langle \sigma \rangle}^W - \varrho)
\]
as $\mathbb {C} W$-modules, where $X_W$ is a set of representatives of
$W$-conjugacy classes of so called special involutions in $W$, $\varrho$ is 
the regular representation of $W$, and $1_{\langle \sigma \rangle}^W$ is the
$\mathbb {C} W$-module induced from the trivial 
$\mathbb {C} \langle \sigma \rangle$-module.
This formula can be deduced 
from earlier work of Lehrer~\cite{Lehrer1987a,Lehrer1987b} and
Fleischmann-Janiszczak~\cite{FleischmannJaniszczak1991,%
FleischmannJaniszczak1993}.
%, see the references in \cite{FelderVeselov2004}.
The main contribution in \cite{FelderVeselov2004} to the theory is a 
uniform geometric description of the sets $X_W$ of 
$W$-conjugacy classes of special involutions used in the formula above.

Felder and Veselov give a similar formula for the twisted action 
where the summation is taken over the set of even elements from $X_W$.

In this note we give a short intrinsic characterisation 
of special involutions in terms of bulky parabolic subgroups.

\section{Notation and Preliminaries}
\label{sec:notation}

Throughout, $W$ denotes a finite  Coxeter group, generated by a set of
simple  reflections   $S  \subseteq  W$;   see~\cite{Bourbaki1968}  or
\cite{GePf2000} for a general  introduction into the theory of Coxeter
groups.  For $J  \subseteq S$, let $W_J$ be  the parabolic subgroup of
$W$ generated  by $J$ and  $w_J$ denotes the  unique word in  $W_J$ of
maximal length (with respect to $S$).  Let $T = S^W$ be the set of all
reflections of  $W$.  Let $\Phi$ be  a root system  with Coxeter group
$W$  and $\Phi_J$  is the  root subsystem  of $\Phi$  corresponding to
$W_J$.  Set  $V :=  {\Z}\Phi \otimes_{\Z} \R$.   Then $V$  affords the
usual reflection  representation of $W$.  For  each involution $\sigma
\in W$  we have a  direct sum decomposition  $V = V_1  \oplus V_{-1}$,
where $V_1$  and $V_{-1}$ are the  $1$ and $-1$-eigenspaces  of $V$ of
$\sigma$, respectively.  For $\epsilon  = \pm 1$ let $\Phi_\epsilon :=
\Phi \cap V_\epsilon$.   Note that for $s = w_J$  we have $\Phi_{-1} =
\Phi_J$.    Following  \cite{FelderVeselov2004},   we   say  that   an
involution $\sigma$ in $W$ is  \emph{special}, if for any root $\alpha
\in  \Phi$  at least  one  of  its  projections onto  $V_\epsilon$  is
proportional to  a root in $\Phi_\epsilon$.   Clearly, this definition
does not depend on the choice of root system for $W$.

The conjugacy  classes of involutions  in $W$ have been  classified by
Richardson~\cite[Thm.\  A]{Richardson1982} and
Springer~\cite{Springer1982} in  terms of  the parabolic
subgroups of  $W$ whose longest  element is central.   More precisely,
each  involution is  conjugate to  a  longest element  $w_J$ which  is
central in $W_J$.

 The normalisers of
parabolic subgroups  of finite Coxeter  groups have been  described by
Howlett~\cite{Howlett1980}        and              Brink        and
Howlett~\cite{BrinkHowlett1999}.  Accordingly, the normaliser $N_W(W_J)$
of $W_J$  in $W$  is a  semi-direct product of  the form  $W_J \rtimes
N_J$, where $N_J$  is itself a semi-direct product  of a Coxeter group
of known type  and a group $M_J$, \cite[Cor.\ 7]{Howlett1980}. 
 It turns out,  however, that in the
case where $w_J$ is central in $W_J$ this group $M_J$ is trivial.

\begin{proposition} \label{p:1}
  $M_J$  acts faithfully as  inner graph  automorphisms on  $W_J$.  In
  particular, if $w_J$ is central in $W_J$, then  $M_J = \{1\}$.
\end{proposition}

\begin{proof}
  According  to the tables  in \cite{BrinkHowlett1999},  a non-trivial
  generator of the group $M_J$ arises from a situation, where $\Size{S
    \setminus J} = 2$ and either $W$  is of type $E_7$ and $W_J$ is of
  type $A_4  \times A_1$, or $W$ is  of type $D_{2n}$ and  $W_J$ is of
  type $A_{2n-2}$ or $A_{2k} \times A_{2l}$ with $k \neq l$ and $k + l
  = n -  1$.  Let us say that $W_J$  is an \emph{M-parabolic} subgroup
  of $W$  in such a case.   An easy check  shows that, if $W_J$  is an
  M-parabolic subgroup of $W$  then $M_J$ induces the same non-trivial
  graph automorphism on $W_J$ as conjugation by $w_J$.  In general, it
  follows that $M_J$ is trivial unless  a conjugate $L$ of $J$ lies in
  a subset  $K \subseteq S$ such  that $\Size{K \setminus L}  = 2$ and
  $W_K = W_N \times W_{K'}$ and $W_L = W_N \times W_{L'}$ for suitable
  subsets  $K', L',  N \subseteq  S$  and $W_{L'}$  is an  M-parabolic
  subgroup  of  $W_{K'}$.   Now  $M_L$  induces  a  non-trivial  inner
  automorphism on $W_L$ and so does $M_J$ on $W_J$.

  By~\cite[Cor.\ 9]{Howlett1980}, $M_J$ intersects the centraliser
  of $J$ in $N_J$ trivially, and hence acts faithfully on $W_J$.
\end{proof}

The  centraliser of  the involution  $w_J$ and  the normaliser  of the
parabolic  subgroup  $W_J$  of  $W$  coincide;  see  
\cite[Prop.\ 7]{FelderVeselov2004}.  We give a new proof of this property.

\begin{proposition}
  Let $J \subseteq S$ be such that $w_J$ is central in $W_J$.  Then
  $C_W(w_J) = N_W(W_J)$.
\end{proposition}

\begin{proof}
  Clearly,  $W_J \subseteq  C_W(w_J)  \cap N_W(W_J)$.   It  thus suffices  to
  consider the set $D_J = \{x \in W  : l(sx) > l(x) \text{ and } l(xs) > l(x)
  \text{ for all  } s \in J\}$ of  distinguished double coset representatives
  of $W_J$ in $W$.
  
  We have $l(w^x) = l(w)$ for all $w \in
  W_J$, $x  \in N_J = \{x \in D_J : J^x = J\}$.  In particular, $w_J^x  = w_J$ for $x  \in N_J$.  Hence
  $N_W(W_J) \subseteq C_W(w_J)$.
  
  Conversely, let $x \in C_W(w_J) \cap  D_J$.  Then $w_J \in W_J \cap W_J^x =
  W_{J \cap J^x}$; cf.~\cite[(2.1.12)]{GePf2000}.   It follows that $J = J^x$
  whence $C_W(w_J) \subseteq N_W(W_J)$.
\end{proof}

We  call  the  parabolic  subgroup  $W_J$  \emph{bulky}  (in  $W$)  if
$N_W(W_J) =  W_J \times  N_J$, i.e., if $N_J$ acts trivially  on $W_J$. 
The  main result of  this note  is the following theorem.

\begin{theorem}
\label{thm:1}
Let $J \subseteq S$ be such  that $w_J$ is central in $W_J$.  Then the
involution $w_J$ is special if and only if $W_J$ is bulky.
\end{theorem}

In  our  arguments  we  do  make  use of  the  classification  of  the
irreducible Coxeter  groups and the  structure of the root  systems of
Weyl groups.   Also, we use the  notation and labelling  of the Dynkin
diagram of $W$ as in \cite[Planches I - IX]{Bourbaki1968}.

%%%%%%%%%%%%%%%%%%%%%%%%%%%%%%%%%%%%%%%%%%%%%%%%%%%%%%%%%%%%%%%%%%%%%%%%%%%%%
\section{Special Involutions and Bulky Parabolic Subgroups}
\label{sec:main}

We maintain the notation from the previous sections.

\begin{lemma}\label{la:1}
  If $\dim V_1 = 1$ and  $\Phi_1 \neq \emptyset$ or if $\dim V_{-1} =
  1$,  then $w_J$  is  special.  In  particular,  
$\pm  s$ is special for every  reflection $s \in T$.
\end{lemma}

\begin{proof}
  The projection of any root onto a one-dimensional space generated by
  a root $\alpha$ is clearly proportional to $\alpha$.
\end{proof}

\begin{rem}\label{rem:autos}
  The element  $w_J$ is central in $W_J$  if and only if  $W_J$ has no
  components of type $A_n$ with $n \geq 2$, of type $D_{2n+1}$ with $n
  \geq 2$, of type $E_6$, or of type $I_2(2m+1)$, $m \geq 2$; see~\cite[1.12]{Richardson1982}.
\end{rem}

\begin{proof}[Proof of Theorem~\ref{thm:1}.]
  We may  assume that  $W$ is irreducible. 
%  $\underline{\Rightarrow}$:
  By \cite[Thm.\ B]{BrinkHowlett1999} and our Proposition~\ref{p:1}
  the group $N_J$ is generated by certain conjugates of elements of
  the form $w_L w_K$, where $L \subseteq K \subseteq S$ such that $L$
  is a conjugate of $J$, $\Size{K \setminus L} = 1$ and $L^{w_K} = L$.
  If $s^{w_L w_K} = s$ for all $s \in L$, then $w_L w_K$ centralises
  $W_L$ and so its conjugate centralises $W_J$.  Obviously, $s^{w_L
    w_K} = s^{w_K}$ for all $s \in L$, since $w_L$ is central in $W_L$.
  
  Now suppose that $W_J$ is not bulky in $W$, i.e., that $N_J$ does not
  centralise $W_J$.   Then there exists a  conjugate $L$ of  $J$ and a
  subset  $K \subseteq  S$ such  that  $L \subseteq  K$ with  $\Size{K
    \setminus  L}   =  1$  and  $w_K$  induces   a  non-trivial  graph
  automorphism on $W_L$.
  It  follows that  $W_K  = W_N  \times  W_{K'}$ for  suitable $N,  K'
  \subseteq S$  where the type of  $W_{K'}$ is one of  those listed in
  Remark~\ref{rem:autos}.
  Since $w_L$ is  central in $W_L$, it follows that  $W_L = W_N \times
  W_{L'}$  where $W_{L'}$  is a  product  of components  of types  not
  listed  in   Remark~\ref{rem:autos}.   Inspection  of   the  maximal
  parabolic subgroups of  $W_{K'}$ shows that $W_{L'}$ is of  type
  $D_{2n}$ and  $W_{K'}$  is of  type  $D_{2n+1}$, $n  \geq 1$;  this
  includes the case where $A_1^2$ embeds into $A_3$ for $n = 1$.
  
  Without loss we may assume $N = \emptyset$, $K' = K = S$ and $L' = L
  = J$.  So let 
\[
  \Phi = \{\pm \varepsilon_i \pm \varepsilon_j : 1 \leq i < j \leq 2n+1\}
\]
 be  a root system of type $D_{2n+1}$ and consider
  the  simple root $\alpha  = \varepsilon_1  - \varepsilon_2$.   It is
  easy to check that $\{\varepsilon_2, \dots, \varepsilon_{2n+1}\}$ is
  a basis  of $V_{-1}$ and  $V_1$ is the $\R$-span  of $\varepsilon_1$
  and  obviously $V_1$ contains no  root, hence  $\Phi_1 =  \Phi \cap  V_1 =
  \emptyset$.   Consequently, no  projection of  a root  in  $\Phi$ on
  $V_1$ is proportional to a root  in $\Phi_1$.  Next we show that the
  projection of $\alpha$ onto $V_{-1}$ is not proportional to any root
  in  $\Phi_{-1}  =  \Phi  \cap  V_{-1}$.  Recall  that  $\Phi_{-1}  =
  \Phi_J$, a root system of type $D_{2n}$ consisting of the roots $\pm
  \varepsilon_i \pm  \varepsilon_j$, with $2  \leq i <  j \leq  2n+1$.  The
  projection of  $v \in  V$ onto  $V_{-1}$ is given  by $\frac12  (v -
  w_J(v))$.    Hence   $\alpha$   projects  onto   $\frac12(\alpha   -
  w_J(\alpha)) = -\varepsilon_2$ and  therefore is not proportional to
  any root in $\Phi_{-1}$.  Thus $w_J$ is not special, as required.

%$\underline{\Leftarrow}$:
  For the converse, suppose now that  $W_J$ is bulky in $W$, i.e., that
  $N_J$ acts trivially  on $W_J$.  If $J = \emptyset$ or  $J = S$, then
  clearly $w_J$  is a  special involution.  So  let us assume  that $J
  \neq \emptyset, S$.  We  consider the different types of irreducible
  Coxeter groups in turn.
  
  If   $W$    is   of   type    $A_n$   ($n   \geq   1$),    then,   by
  Remark~\ref{rem:autos},  $W_J$ is  necessarily a  direct  product of
  components  of type  $A_1$.   But if  there  is more  than one  such
  component,  $N_J$ permutes  them non-trivially.   Hence $W_J$  is of
  type $A_1$ and the claim follows by Lemma~\ref{la:1}.
  
  If   $W$    is   of   type    $C_n$   ($n   \geq   2$),    then,   by
  Remark~\ref{rem:autos}, $W_J$ is a  direct product of a component of
  type $C_m$, $0 \leq m <  n$ and further components of type $A_1$. 
  As before  there cannot be  more than one  component of type  $A_1$. 
  Hence $W_J$ is of type $C_m$ or of type $C_m \times A_1$ for some $m
  < n$.  In any case $W_J$ has a component of type $C_m$.

Let  
\[
  \Phi  = \{\pm2\varepsilon_i  :  1 \leq  i  \leq  n\} \cup  \{\pm
\varepsilon_i \pm  \varepsilon_j : 1 \leq  i < j \leq n\}
\]
be the root
system of type $C_n$.  Consider  the maximal rank subsystem $\Phi'$ of
type $C_m  \times C_{n-m}$ consisting of  the long roots
$\{\pm2\varepsilon_i : 1  \leq i \leq n\}$ and  the short roots $\{\pm
\varepsilon_i \pm \varepsilon_j  : 1 \leq i < j \leq m  \text{ or } m+1
\leq  i <  j  \leq n\}$.  Let  $U_1$  be  the subspace  of $V$ spanned  by
$\varepsilon_1, \dots,  \varepsilon_m$ and $U_2$  the subspace spanned
by  $\varepsilon_{m+1}, \dots, \varepsilon_n$.  Then $U_2 \cap \Phi$
is a root system of type $C_{n-m}$.
All the  long roots  $\pm 2\varepsilon_i$ of  $\Phi$ are  contained in
$\Phi'$.  A short root $\pm \varepsilon_i \pm \varepsilon_j$ is either
contained in  $\Phi'$ or both its  projections on $U_1$  and $U_2$ are
proportional to a root in $\Phi'$.
By construction, the $-1$-eigenspace $V_{-1}$ of $w_J$ contains $U_1$.
Hence every  root that lies in $U_1$  or is proportional to  a root in
$U_1$  is also  proportional to  a root  in $V_{-1}$.   It  remains to
consider the roots in $U_2$.  Without loss we can now assume that $m =
0$.   Then  $W_J$   is  of  type  $A_1$  and   the  claim  follows  by
Lemma~\ref{la:1}.

If   $W$   is   of   type   $D_{2n+1}$   ($n   \geq   2$),   then,   by
Remark~\ref{rem:autos},  $W_J$  is a  direct  product  of an  optional
component of type  $D_{2m}$, $1 \leq m \leq  n$ and further components
of type $A_1$.   As before there cannot be more  than one component of
type  $A_1$.  And  $N_J$ acts  non-trivially  on a  component of  type
$D_{2m}$.   Hence $W_J$  is of  type $A_1$  and the  claim  follows by
Lemma~\ref{la:1}.

If   $W$   is    of   type   $D_{2n}$   ($n   \geq    2$),   then,   by
Remark~\ref{rem:autos},  $W_J$  is a  direct  product  of an  optional
component of type  $D_{2m}$, $1 \leq m <  n$ and further components
of type $A_1$.   As before there cannot be more  than one component of
type $A_1$.  The non-trivial action  of the parabolic subgroup of type
$D_{2m+1}$ on a component of  type $D_{2m}$ then restricts the type of
$W_J$ to either $A_1$ or  $D_{2(n-1)} \times A_1$.  In the latter case
$V_1 \cap \Phi$  is a root system  of type $A_1$ and so  in both cases
the claim follows by Lemma~\ref{la:1}.

If $W$ is  of type $I_2(m)$ ($m  \geq 5$), then $W_J$ is  of type $A_1$
and the claim follows by Lemma~\ref{la:1}.

Finally, if  $W$ has type  $E_6$, $E_7$,
$E_8$,  $F_4$,  $H_3$, or  $H_4$,  then  the  claim is  established  by
inspection.
\end{proof}

\begin{rem}
  Felder and Veselov prove  the one implication of Theorem \ref{thm:1},
  namely that $W_J$ is  bulky if $w_J$ is special in  a case  by case
  analysis~\cite[Prop.\ 10]{FelderVeselov2004}.
\end{rem}

\begin{rem}
  Bulky parabolic  subgroups can be  easily classified.  It  turns out
  that if $W$ has a central  longest element, then $w_J$ is central in
  $W_J$ whenever  $W_J$ is bulky.  Otherwise, $W$  has bulky parabolic
  subgroups  $W_J$  which  are  not  associated  with  a  conjugacy  
class of involutions in $W$.
\end{rem}

%%%%%%%%%%%%%%%%%%%%%%%%%%%%%%%%%%%%%%%%%%%%%%%%%%%%%%%%%%%%%%%%%%%%%%%%%%%%%
\bigskip

\noindent
{\bf Acknowledgements}: This paper  was written while the first author
was  visiting  the  School   of  Mathematics  and  Statistics  of  the
University of Birmingham under the Scheme 2 LMS grant no.~2925. We are
grateful to the LMS for its financial support and to the members of the
School for their  hospitality.  We also wish to  thank A.\ Borovik for
bringing the  problem  of  an intrinsic  characterisation  of  special
involutions to our attention.

%%%%%%%%%%%%%%%%%%%%%%%%%%%%%%%%%%%%%%%%%%%%%%%%%%%%%%%%%%%%%%%%%%%%%%%%%%%%%
\bibliographystyle{amsplain}
\bibliography{descent}

\end{document}